\documentclass{amsproc}
\usepackage{mathtools}
\usepackage{amsmath}
\usepackage[dvipsnames]{xcolor}
\usepackage{tikz}
\usepackage{listings}
\usepackage{courier}
\usepackage{siunitx}
\usepackage{ifthen}
\newif\ifcuboidshade
\newif\ifcuboidemphedge

\newcommand{\Xsub}[1]{\mathbin{\mathop{\times}\limits_{#1}}}

\usepackage{anyfontsize}

\definecolor{Y}{RGB}{255,255,110}
\definecolor{R}{RGB}{255,70,75}
\definecolor{G}{RGB}{151,216,56}
\definecolor{B}{RGB}{51,152,237}
\definecolor{W}{RGB}{255,255,255}
\definecolor{O}{RGB}{255,165,0}

\newcommand{\gridthick}{0.6pt}
\newcommand{\borderthick}{1.2pt}
\newcommand{\numsize}{\small\ttfamily} 
\newcommand{\centerfillopacity}{0.35}  
\newif\ifshowfacelabel \showfacelabelfalse 

\definecolor{FaceU}{named}{G}
\definecolor{FaceD}{named}{Y}
\definecolor{FaceF}{named}{B}
\definecolor{FaceB}{named}{W}
\definecolor{FaceR}{named}{R}
\definecolor{FaceL}{named}{O}

\newcommand{\drawFaceRows}[9]{%
  \begin{scope}[shift={(#1,#2)}]
    \draw[line width=\borderthick] (0,0) rectangle (5,5);
    \foreach \i in {1,...,4} {%
      \draw[line width=\gridthick] (\i,0) -- (\i,5);
      \draw[line width=\gridthick] (0,\i) -- (5,\i);
    }
    \fill[fill=#4, fill opacity=\centerfillopacity, draw=none] (0,0) rectangle (5,5);
    
    \ifshowfacelabel
      \node[font=\bfseries] at (2.5, 5.45) {#3};
    \fi

    \foreach \t [count=\c] in {#5} { \node at (\c-0.5, 4.5) {\numsize \t}; }
    \foreach \t [count=\c] in {#6} { \node at (\c-0.5, 3.5) {\numsize \t}; }
    \foreach \t [count=\c] in {#7} { \node at (\c-0.5, 2.5) {\numsize \t}; }
    \foreach \t [count=\c] in {#8} { \node at (\c-0.5, 1.5) {\numsize \t}; }
    \foreach \t [count=\c] in {#9} { \node at (\c-0.5, 0.5) {\numsize \t}; }
  \end{scope}
}

\tikzset{
  cuboid/.is family,
  cuboid,
  shiftx/.initial=0,
  shifty/.initial=0,
  dimx/.initial=3,
  dimy/.initial=3,
  dimz/.initial=3,
  scale/.initial=1,
  densityx/.initial=1,
  densityy/.initial=1,
  densityz/.initial=1,
  rotation/.initial=0,
  anglex/.initial=0,
  angley/.initial=90,
  anglez/.initial=225,
  scalex/.initial=1,
  scaley/.initial=1,
  scalez/.initial=0.5,
  front/.style={draw=black,fill=B},
  top/.style={draw=black,fill=G},
  right/.style={draw=black,fill=R},
  shade/.is if=cuboidshade,
  shadecolordark/.initial=black,
  shadecolorlight/.initial=white,
  shadeopacity/.initial=0.15,
  shadesamples/.initial=16,
  emphedge/.is if=cuboidemphedge,
  emphstyle/.style={thick},
}

\newcommand{\colorcubiefront}[4][]{%
  \begingroup
  \tikzset{cuboid,#1}%
  \pgfmathsetlengthmacro{\vectorxx}{\tikzcuboidkey{scalex}*cos(\tikzcuboidkey{anglex})*28.452756}
  \pgfmathsetlengthmacro{\vectorxy}{\tikzcuboidkey{scalex}*sin(\tikzcuboidkey{anglex})*28.452756}
  \pgfmathsetlengthmacro{\vectoryx}{\tikzcuboidkey{scaley}*cos(\tikzcuboidkey{angley})*28.452756}
  \pgfmathsetlengthmacro{\vectoryy}{\tikzcuboidkey{scaley}*sin(\tikzcuboidkey{angley})*28.452756}
  \pgfmathsetlengthmacro{\vectorzx}{\tikzcuboidkey{scalez}*cos(\tikzcuboidkey{anglez})*28.452756}
  \pgfmathsetlengthmacro{\vectorzy}{\tikzcuboidkey{scalez}*sin(\tikzcuboidkey{anglez})*28.452756}
  \begin{scope}[
    xshift=\tikzcuboidkey{shiftx}, yshift=\tikzcuboidkey{shifty},
    scale=\tikzcuboidkey{scale}, rotate=\tikzcuboidkey{rotation},
    x={(\vectorxx,\vectorxy)}, y={(\vectoryx,\vectoryy)}, z={(\vectorzx,\vectorzy)}
  ]
    \pgfmathsetmacro{\sx}{1/\tikzcuboidkey{densityx}}
    \pgfmathsetmacro{\sy}{1/\tikzcuboidkey{densityy}}
    \pgfmathsetmacro{\ix}{#2}
    \pgfmathsetmacro{\iy}{#3}
    \pgfmathsetmacro{\lowx}{(\ix-1)*\sx}
    \pgfmathsetmacro{\x}{\ix*\sx}
    \pgfmathsetmacro{\lowy}{(\iy-1)*\sy}
    \pgfmathsetmacro{\y}{\iy*\sy}
    \fill[draw=black, line width=0.4pt, fill=#4]
      (\lowx,\lowy,\tikzcuboidkey{dimz}) --
      (\lowx,\y,\tikzcuboidkey{dimz}) --
      (\x,\y,\tikzcuboidkey{dimz}) --
      (\x,\lowy,\tikzcuboidkey{dimz}) -- cycle;
  \end{scope}
  \endgroup
}

\newcommand{\colorcubieright}[4][]{%
  \begingroup
  \tikzset{cuboid,#1}%
  \pgfmathsetlengthmacro{\vectorxx}{\tikzcuboidkey{scalex}*cos(\tikzcuboidkey{anglex})*28.452756}
  \pgfmathsetlengthmacro{\vectorxy}{\tikzcuboidkey{scalex}*sin(\tikzcuboidkey{anglex})*28.452756}
  \pgfmathsetlengthmacro{\vectoryx}{\tikzcuboidkey{scaley}*cos(\tikzcuboidkey{angley})*28.452756}
  \pgfmathsetlengthmacro{\vectoryy}{\tikzcuboidkey{scaley}*sin(\tikzcuboidkey{angley})*28.452756}
  \pgfmathsetlengthmacro{\vectorzx}{\tikzcuboidkey{scalez}*cos(\tikzcuboidkey{anglez})*28.452756}
  \pgfmathsetlengthmacro{\vectorzy}{\tikzcuboidkey{scalez}*sin(\tikzcuboidkey{anglez})*28.452756}
  \begin{scope}[
    xshift=\tikzcuboidkey{shiftx}, yshift=\tikzcuboidkey{shifty},
    scale=\tikzcuboidkey{scale}, rotate=\tikzcuboidkey{rotation},
    x={(\vectorxx,\vectorxy)}, y={(\vectoryx,\vectoryy)}, z={(\vectorzx,\vectorzy)}
  ]
    \pgfmathsetmacro{\sy}{1/\tikzcuboidkey{densityy}}
    \pgfmathsetmacro{\sz}{1/\tikzcuboidkey{densityz}}
    \pgfmathsetmacro{\jy}{#2}
    \pgfmathsetmacro{\kz}{#3}
    \pgfmathsetmacro{\lowy}{(\jy-1)*\sy}
    \pgfmathsetmacro{\y}{\jy*\sy}
    \pgfmathsetmacro{\lowz}{(\kz-1)*\sz}
    \pgfmathsetmacro{\z}{\kz*\sz}
    \fill[draw=black, line width=0.4pt, fill=#4]
      (\tikzcuboidkey{dimx},\lowy,\lowz) --
      (\tikzcuboidkey{dimx},\lowy,\z) --
      (\tikzcuboidkey{dimx},\y,\z) --
      (\tikzcuboidkey{dimx},\y,\lowz) -- cycle;
  \end{scope}
  \endgroup
}

\newcommand{\colorcubietop}[4][]{%
  \begingroup
  \tikzset{cuboid,#1}%
  \pgfmathsetlengthmacro{\vectorxx}{\tikzcuboidkey{scalex}*cos(\tikzcuboidkey{anglex})*28.452756}
  \pgfmathsetlengthmacro{\vectorxy}{\tikzcuboidkey{scalex}*sin(\tikzcuboidkey{anglex})*28.452756}
  \pgfmathsetlengthmacro{\vectoryx}{\tikzcuboidkey{scaley}*cos(\tikzcuboidkey{angley})*28.452756}
  \pgfmathsetlengthmacro{\vectoryy}{\tikzcuboidkey{scaley}*sin(\tikzcuboidkey{angley})*28.452756}
  \pgfmathsetlengthmacro{\vectorzx}{\tikzcuboidkey{scalez}*cos(\tikzcuboidkey{anglez})*28.452756}
  \pgfmathsetlengthmacro{\vectorzy}{\tikzcuboidkey{scalez}*sin(\tikzcuboidkey{anglez})*28.452756}
  \begin{scope}[
    xshift=\tikzcuboidkey{shiftx}, yshift=\tikzcuboidkey{shifty},
    scale=\tikzcuboidkey{scale}, rotate=\tikzcuboidkey{rotation},
    x={(\vectorxx,\vectorxy)}, y={(\vectoryx,\vectoryy)}, z={(\vectorzx,\vectorzy)}
  ]
    \pgfmathsetmacro{\sx}{1/\tikzcuboidkey{densityx}}
    \pgfmathsetmacro{\sz}{1/\tikzcuboidkey{densityz}}
    \pgfmathsetmacro{\ix}{#2}
    \pgfmathsetmacro{\kz}{#3}
    \pgfmathsetmacro{\lowx}{(\ix-1)*\sx}
    \pgfmathsetmacro{\x}{\ix*\sx}
    \pgfmathsetmacro{\lowz}{(\kz-1)*\sz}
    \pgfmathsetmacro{\z}{\kz*\sz}
    \fill[draw=black, line width=0.4pt, fill=#4]
      (\lowx,\tikzcuboidkey{dimy},\lowz) --
      (\lowx,\tikzcuboidkey{dimy},\z) --
      (\x,\tikzcuboidkey{dimy},\z) --
      (\x,\tikzcuboidkey{dimy},\lowz) -- cycle;
  \end{scope}
  \endgroup
}

\newcommand{\tikzcuboidkey}[1]{\pgfkeysvalueof{/tikz/cuboid/#1}}

\newcommand{\tikzcuboid}[1]{
    \tikzset{cuboid,#1} 
  \pgfmathsetlengthmacro{\vectorxx}{\tikzcuboidkey{scalex}*cos(\tikzcuboidkey{anglex})*28.452756}
  \pgfmathsetlengthmacro{\vectorxy}{\tikzcuboidkey{scalex}*sin(\tikzcuboidkey{anglex})*28.452756}
  \pgfmathsetlengthmacro{\vectoryx}{\tikzcuboidkey{scaley}*cos(\tikzcuboidkey{angley})*28.452756}
  \pgfmathsetlengthmacro{\vectoryy}{\tikzcuboidkey{scaley}*sin(\tikzcuboidkey{angley})*28.452756}
  \pgfmathsetlengthmacro{\vectorzx}{\tikzcuboidkey{scalez}*cos(\tikzcuboidkey{anglez})*28.452756}
  \pgfmathsetlengthmacro{\vectorzy}{\tikzcuboidkey{scalez}*sin(\tikzcuboidkey{anglez})*28.452756}
  \begin{scope}[xshift=\tikzcuboidkey{shiftx}, yshift=\tikzcuboidkey{shifty}, scale=\tikzcuboidkey{scale}, rotate=\tikzcuboidkey{rotation}, x={(\vectorxx,\vectorxy)}, y={(\vectoryx,\vectoryy)}, z={(\vectorzx,\vectorzy)}]
    \pgfmathsetmacro{\steppingx}{1/\tikzcuboidkey{densityx}}
  \pgfmathsetmacro{\steppingy}{1/\tikzcuboidkey{densityy}}
  \pgfmathsetmacro{\steppingz}{1/\tikzcuboidkey{densityz}}
  \newcommand{\dimx}{\tikzcuboidkey{dimx}}
  \newcommand{\dimy}{\tikzcuboidkey{dimy}}
  \newcommand{\dimz}{\tikzcuboidkey{dimz}}
  \pgfmathsetmacro{\secondx}{2*\steppingx}
  \pgfmathsetmacro{\secondy}{2*\steppingy}
  \pgfmathsetmacro{\secondz}{2*\steppingz}
  \ifthenelse{\equal{\dimx}{1}}
    {\foreach \x in {\steppingx,...,\dimx}}
    {\foreach \x in {\steppingx,\secondx,...,\dimx}}
  {     \ifthenelse{\equal{\dimy}{1}}
    {\foreach \y in {\steppingy,...,\dimy}}
    {\foreach \y in {\steppingy,\secondy,...,\dimy}}
    {   \pgfmathsetmacro{\lowx}{(\x-\steppingx)}
      \pgfmathsetmacro{\lowy}{(\y-\steppingy)}
      \filldraw[cuboid/front] (\lowx,\lowy,\dimz) -- (\lowx,\y,\dimz) -- (\x,\y,\dimz) -- (\x,\lowy,\dimz) -- cycle;
    }
    }
    \ifthenelse{\equal{\dimx}{1}}
    {\foreach \x in {\steppingx,...,\dimx}}
    {\foreach \x in {\steppingx,\secondx,...,\dimx}}
  { \ifthenelse{\equal{\dimz}{1}}
    {\foreach \z in {\steppingz,...,\dimz}}
    {\foreach \z in {\steppingz,\secondz,...,\dimz}}
    {   \pgfmathsetmacro{\lowx}{(\x-\steppingx)}
      \pgfmathsetmacro{\lowz}{(\z-\steppingz)}
      \filldraw[cuboid/top] (\lowx,\dimy,\lowz) -- (\lowx,\dimy,\z) -- (\x,\dimy,\z) -- (\x,\dimy,\lowz) -- cycle;
        }
    }
    \ifthenelse{\equal{\dimy}{1}}
    {\foreach \y in {\steppingy,...,\dimy}}
    {\foreach \y in {\steppingy,\secondy,...,\dimy}}
  { \ifthenelse{\equal{\dimz}{1}}
    {\foreach \z in {\steppingz,...,\dimz}}
    {\foreach \z in {\steppingz,\secondz,...,\dimz}}
    {   \pgfmathsetmacro{\lowy}{(\y-\steppingy)}
      \pgfmathsetmacro{\lowz}{(\z-\steppingz)}
      \filldraw[cuboid/right] (\dimx,\lowy,\lowz) -- (\dimx,\lowy,\z) -- (\dimx,\y,\z) -- (\dimx,\y,\lowz) -- cycle;
    }
  }
  \ifcuboidemphedge
    \draw[cuboid/emphstyle] (0,\dimy,0) -- (\dimx,\dimy,0) -- (\dimx,\dimy,\dimz) -- (0,\dimy,\dimz) -- cycle;%
    \draw[cuboid/emphstyle] (0,\dimy,\dimz) -- (0,0,\dimz) -- (\dimx,0,\dimz) -- (\dimx,\dimy,\dimz);%
    \draw[cuboid/emphstyle] (\dimx,\dimy,0) -- (\dimx,0,0) -- (\dimx,0,\dimz);%
    \fi

    \ifcuboidshade
    \pgfmathsetmacro{\cstepx}{\dimx/\tikzcuboidkey{shadesamples}}
    \pgfmathsetmacro{\cstepy}{\dimy/\tikzcuboidkey{shadesamples}}
    \pgfmathsetmacro{\cstepz}{\dimz/\tikzcuboidkey{shadesamples}}
    \foreach \s in {1,...,\tikzcuboidkey{shadesamples}}
    {   \pgfmathsetmacro{\lows}{\s-1}
        \pgfmathsetmacro{\cpercent}{(\lows)/(\tikzcuboidkey{shadesamples}-1)*100}
        \fill[opacity=\tikzcuboidkey{shadeopacity},color=\tikzcuboidkey{shadecolorlight}!\cpercent!\tikzcuboidkey{shadecolordark}] (0,\s*\cstepy,\dimz) -- (\s*\cstepx,\s*\cstepy,\dimz) -- (\s*\cstepx,0,\dimz) -- (\lows*\cstepx,0,\dimz) -- (\lows*\cstepx,\lows*\cstepy,\dimz) -- (0,\lows*\cstepy,\dimz) -- cycle;
        \fill[opacity=\tikzcuboidkey{shadeopacity},color=\tikzcuboidkey{shadecolorlight}!\cpercent!\tikzcuboidkey{shadecolordark}] (0,\dimy,\s*\cstepz) -- (\s*\cstepx,\dimy,\s*\cstepz) -- (\s*\cstepx,\dimy,0) -- (\lows*\cstepx,\dimy,0) -- (\lows*\cstepx,\dimy,\lows*\cstepz) -- (0,\dimy,\lows*\cstepz) -- cycle;
        \fill[opacity=\tikzcuboidkey{shadeopacity},color=\tikzcuboidkey{shadecolorlight}!\cpercent!\tikzcuboidkey{shadecolordark}] (\dimx,0,\s*\cstepz) -- (\dimx,\s*\cstepy,\s*\cstepz) -- (\dimx,\s*\cstepy,0) -- (\dimx,\lows*\cstepy,0) -- (\dimx,\lows*\cstepy,\lows*\cstepz) -- (\dimx,0,\lows*\cstepz) -- cycle;
    }
    \fi 

  \end{scope}
}

\makeatother

\newcommand{\pare}[1]{\left( #1 \right)} 

\newcommand{\sign}{\operatorname{sign}}

\newtheorem{theorem}{Theorem}
\newtheorem{question}{Question}
\theoremstyle{definition}

\newtheorem{remark}{Remark}

\def\Sym{\mathbb{S}}

\def\A{\mathbb{A}}

\def\Q{\mathbb{Q}}

\def\Z{\mathbb{Z}}
\def\S{\mathbb{S}}

\def\id{\operatorname{id}}

\def\Gal{{\operatorname{Gal}}}
\def\sgn{{\operatorname{sign}\,}}

\definecolor{Y}{RGB}{255,255,110}
\definecolor{R}{RGB}{255,70,75}
\definecolor{G}{RGB}{151,216,56}
\definecolor{B}{RGB}{51,152,237}
\definecolor{W}{RGB}{255,255,255}
\definecolor{O}{RGB}{255,165,0}

\newcommand{\Cyc}[1]{ \Z/(#1) }
\newcommand{\Rub}{\mathcal{R}}

\newcommand{\Galf}[2]{ \Gal\left(#1,#2\right)}

\newcommand{\disc}[1]{\operatorname{disc}({#1})}

\begin{document}

\title{Galois' Professor's Revenge}

\markright{Submission}


\author{M. Damele}
\author{A. Loi}
\author{M. Mereb}
\author{L. Vendramin}

\address[Damele and Loi]{University of Cagliari (Università degli Studi di Cagliari)}
\email{m.damele4@studenti.unica.it}
\email{loi@unica.it}

\address[Mereb]{IMAS--CONICET and Depto. de Matem\'atica, FCEN, Universidad de Buenos Aires, Pab.~1, Ciudad Universitaria, C1428EGA, Buenos Aires, Argentina}
\email{mmereb@dm.uba.ar}

\address[Vendramin]{Department of Mathematics and Data
Science, Vrije Universiteit Brussel, Pleinlaan 2, 1050 Brussels, Belgium}
\email{Leandro.Vendramin@vub.be}

\begin{abstract}
We prove that the groups associated with the Revenge Cube and the Professor's Cube can be realized as Galois groups over the rationals.
\end{abstract}

\maketitle

\section{Introduction}

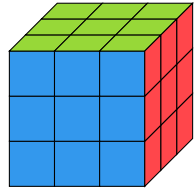
\begin{figure}[h]
\centering 
   \begin{tikzpicture}[scale=0.6]
       \tikzset{
    cuboid/front/.style={draw=black, fill=B},
    cuboid/right/.style={draw=black, fill=R},
    cuboid/top/.style={draw=black, fill=G}}
    \tikzcuboid{dimx=3,dimy=3,dimz=3,shiftx=50,shifty=40};
     \end{tikzpicture}
    \caption{The Rubik's Cube.}
    \label{fig:Rubik}
\end{figure}

In the 1970s, the celebrated Rubik’s Cube was invented (see Figure~\ref{fig:Rubik}). The set of moves of the Rubik’s Cube form a group, which we denote by $\Rub_3$. 
The order of $\Rub_3$ is
\[
43252003274489856000=2^{27}3^{14}5^37^211.
\]
As is well known, there are several mathematical problems related to this group, and numerous excellent references discuss them; see for example~\cite{bandelow2012inside}. 

Less well known, but equally interesting, is the existence of other similar puzzles, including the Revenge Cube (see Figure~\ref{fig:revenge}) and the Professor’s Cube (see Figure~\ref{fig:professor}). As with the Rubik’s Cube, the sets of moves of these puzzles also form groups, denoted by $\Rub_4$ and
$\Rub_5$, respectively. 


The groups $\Rub_4$ and $\Rub_5$ are extraordinarily large. For example,
the order of $\Rub_4$~is 
\begin{align*}
16972688908618238933770849245964147960401887232000000000.
\end{align*}
In the case of the Professor’s Cube, the order is even more astronomical:
\begin{align*}
|\Rub_5| &= |\Rub_3|\frac{(24!)^3}{4} \approx 2.58\times 10^{90}.
\end{align*}
To put this in perspective, the exact value is
\num[group-separator={\ \linebreak[1]}]{2582636272886959379162819698174683585918088940054237132
144778804568925405184000000000000000}.   



These groups, of course, can be used to analyze and solve the puzzles. Moreover, there remain several open problems related to these puzzles---most notably, determining the minimal number of moves required to solve them from any scrambled position. This number is known as God’s number, and for the Rubik’s Cube $\Rub_3$, it is known to be twenty~\cite{zbMATH06212024}. For higher-order puzzles like the Revenge Cube and the Professor’s Cube, however, this problem remains unsolved.

\begin{figure}[ht] 
  \begin{minipage}[b]{0.5\linewidth}
    \centering
   \begin{tikzpicture}[scale=0.4]
    \tikzset{
    cuboid/front/.style={draw=black, fill=B},
    cuboid/right/.style={draw=black, fill=R},
    cuboid/top/.style={draw=black, fill=G}}
   \tikzcuboid{dimx=4,dimy=4,dimz=4,shiftx=50,shifty=40};
    \end{tikzpicture}
    \caption{The Revenge Cube.}
    \label{fig:revenge}
    \vspace{4ex}
  \end{minipage}
  \begin{minipage}[b]{0.5\linewidth}
    \centering
 \begin{tikzpicture}[scale=0.4]
    \tikzset{
    cuboid/front/.style={draw=black, fill=B},
    cuboid/right/.style={draw=black, fill=R},
    cuboid/top/.style={draw=black, fill=G}}
   \tikzcuboid{dimx=5,dimy=5,dimz=5,shiftx=50,shifty=40};
    \end{tikzpicture}
   \caption{The Professor's Cube.}
 \label{fig:professor}
    \vspace{4ex}
  \end{minipage} 
\end{figure}

This paper investigates the realization of these puzzle groups as Galois groups over the rational numbers. This is an instance of the well-known \emph{Inverse Galois Problem}. More concretely, this paper aims to address the following natural question:

\begin{question}
    Let $n\in\{3,4,5\}$. Can the group $\Rub_n$ be realized as a Galois group over the rationals?  
\end{question}

For those unfamiliar with the Inverse Galois Problem, let us briefly explain that it asks whether every finite group can be realized as a Galois group over the field of rational numbers. Several families of groups are known to occur as Galois groups, including all finite solvable groups, symmetric and alternating groups, and various simple groups. However, the problem remains wide open in full generality. For example, it is still unknown whether the Mathieu group $M_{23}$ can be realized as a Galois group over $\Q$. Moreover, the Inverse Galois Problem remains open for most finite simple groups of Lie type. For further information on this fascinating branch of mathematics, we refer the reader to \cite{MR3822366,MR2363329,MR1405612}.

In the case of the classical Rubik’s Cube, the group $\Rub_3$ is known to be realizable as a Galois group over the rationals~\cite{arXiv:2411.11566}. Moreover, it is possible to construct an explicit rational polynomial whose Galois group is isomorphic to $\Rub_3$, as the following theorem shows.

\begin{theorem}
    \label{thm:Rubik}
    Let 
    \[
    g(X) = X^{24} + c^2
    (X^2 + 1)\quad\text{and}\quad 
    f(X)=
 (X^2-X)^8 f_8\left(\frac{ X^3-3X+1 }{X^2-X} \right) 
    \]
where
$f_8(X) = X^8-a X-b .$  
    Then the polynomial
    $fg$ has Galois group over the rationals isomorphic to the Rubik's Cube group $\Rub_3$
    for a suitable specialization of parameters $a,b$ and $c.$
\end{theorem}

In~\cite{arXiv:2411.11566} it is shown that
\begin{equation}
\label{eq:Rubik_parameters}
a= 2139,\quad b=-6489,\quad 
c = 
\frac{1962764241992810496}{619884697145165705}
\end{equation}
is a possible choice to get Rubik's group as Galois' group of the splitting field.


The theorem can be verified using Magma~\cite{zbMATH01077111}; here we use Magma V2.28-18. The computation of the Galois group of $fg$ takes less than two minutes on a desktop computer equipped with an AMD Ryzen 7 7700X 8-core processor and 128 GB of RAM. Magma verifies almost immediately that this Galois group is isomorphic to~$\Rub_3$.

Techniques similar to those used to prove Theorem~\ref{thm:Rubik} can be adapted to 
realize more general Rubik-type groups as Galois 
groups over the rationals. 
Here, we 
present an explicit result for cubes of size four and five. 

\begin{theorem}
\label{thm:n=4}
    Let $f$ be as in Theorem~\ref{thm:Rubik},  
    \begin{align*}
     g=X^{24}-t(X+1)\quad\text{and}\quad 
    h= X^{24}-X-1.
    \end{align*}
    Then 
    $fgh$ has Galois group isomorphic 
    to the Revenge Cube group $\Rub_4$ for a suitable specialization 
    of the parameter $t$. 
\end{theorem}

For instance, in this theorem one can take 
\begin{gather*}
a = 2139,\quad b=-6489,\quad 
t=-\frac{2^{67}3^{24}}{23^{23}\cdot 31\cdot 281\cdot 1201\cdot 70529\cdot 9801219477271}.
\end{gather*}

Theorem~\ref{thm:n=4} can be verified using Magma. The computation of the Galois group takes approximately twelve minutes on our desktop computer. As before, Magma quickly verifies that this Galois group is isomorphic to $\Rub_4$. For the proof of the theorem, including an explanation of the origin of our polynomials, see Section~\ref{section:proofs}.


\begin{theorem}
\label{thm:n=5}
    Let $f$ and $g$ be as in Theorem~\ref{thm:Rubik}. For 
    $i\in\{1,2,3\}$, let 
    \[
    h_i = X^{24} - u_i(X +1). 
    \]
Then $fgh_1h_2h_3$ has Galois group isomorphic 
to the Professor's Cube group $\Rub_5$ for a suitable specialization of the parameters. 
\end{theorem}

For example, in this theorem one can take $a$, $b$ and $c$ as in~\eqref{eq:Rubik_parameters}, and 
\begin{align*}
    u_1 &= -2^{67}\cdot 3^{24}\left(23^{23}\cdot 31\cdot 281\cdot 1201\cdot 70529\cdot 9801219477271\right)^{-1},\\
    u_2 &= 2^{72}\cdot 3^{24}\cdot 7\cdot 1437417619559484462138047\cdot 23^{-22},
\end{align*}
and 
\begin{align*}
        u_3 = \dfrac{2^{75}\cdot 3^{14}\cdot 31\cdot 281\cdot 1201\cdot 70529
        \cdot 9801219477271}{7^{2}\cdot 23^{22}\cdot 195574568093355782014153^2}.
\end{align*}



Magma computed the Galois 
group of $fgh_1h_2h_3$ in approximately three days, while it 
verified in about 40 minutes
that the resulting group is isomorphic to $\Rub_5$, 
using the permutation representation described in Appendix~\ref{section:generators}.

\section{The groups}

Let $f\colon G\to K$ and $g\colon H\to K$ be group homomorphisms. 
The \emph{fiber product} $G\times_{f,g}H$ is the
subgroup 
\[
\{(x,y)\in G\times H:f(x)=g(y)\}
\]
of $G\times H$. 

Let $n,m\geq2$ and $G$ be a subgroup of $\S_m$. Following~\cite{arXiv:2411.11566}, we
write $(\Cyc{n}\wr G)^ \circ$ to denote
the kernel of the group homomorphism
\begin{align*}
\Cyc{n}\wr G \to \Cyc{n}, \quad
(x, \sigma) \mapsto \sum_{i=1}^m x_i,
\end{align*}
where $x=(x_1,\dots,x_m)$.

The structure of the Rubik’s Cube group is well known; see, for example, Joyner’s book~\cite{MR2599606}.
The group $\Rub_3$ is isomorphic to 
  $$
    \pare{\Cyc{3}\wr \S_8}^\circ \times_{\sgn}
    \pare{\Cyc{2}\wr \S_{12}}^\circ,
    $$
    where the symbol $\times_\sgn$ indicates a fiber product with respect to both maps
    $$\Cyc{n}\wr \S_m \to \{\pm 1\},\quad
    (x,\sigma)\mapsto \sgn(\sigma),
    $$
    for $(n,m)\in\{(3,8),(2,12)\}$.

To understand the group structure of $\Rub_3$ 
we refer to the following picture:
\[
\Rub_3\simeq
  \vcenter{\hbox{\begin{tikzpicture}[scale=0.4]
  \tikzset{
    cuboid/front/.style={draw=black, fill=white},
    cuboid/right/.style={draw=black, fill=white},
    cuboid/top/.style={draw=black, fill=white}
  }
  \tikzcuboid{
    dimx=3, dimy=3, dimz=3,
    shiftx=50, shifty=40
  };
  \colorcubiefront[shiftx=50,shifty=40,dimx=3,dimy=3,dimz=3]{3}{3}{B}
  \colorcubiefront[shiftx=50,shifty=40,dimx=3,dimy=3,dimz=3]{1}{3}{B}
  \colorcubiefront[shiftx=50,shifty=40,dimx=3,dimy=3,dimz=3]{3}{1}{B}
  \colorcubiefront[shiftx=50,shifty=40,dimx=3,dimy=3,dimz=3]{1}{1}{B}
  \colorcubietop[shiftx=50,shifty=40,dimx=3,dimy=3,dimz=3]{3}{3}{G}
  \colorcubietop[shiftx=50,shifty=40,dimx=3,dimy=3,dimz=3]{1}{3}{G}
  \colorcubietop[shiftx=50,shifty=40,dimx=3,dimy=3,dimz=3]{3}{1}{G}
  \colorcubietop[shiftx=50,shifty=40,dimx=3,dimy=3,dimz=3]{1}{1}{G}
  \colorcubieright[shiftx=50,shifty=40,dimx=3,dimy=3,dimz=3]{3}{3}{R}
  \colorcubieright[shiftx=50,shifty=40,dimx=3,dimy=3,dimz=3]{1}{3}{R}
  \colorcubieright[shiftx=50,shifty=40,dimx=3,dimy=3,dimz=3]{3}{1}{R}
  \colorcubieright[shiftx=50,shifty=40,dimx=3,dimy=3,dimz=3]{1}{1}{R}
\end{tikzpicture}   }}
  \Xsub{\sign}
  \vcenter{\hbox{\begin{tikzpicture}[scale=0.4]
  \tikzset{
    cuboid/front/.style={draw=black, fill=white},
    cuboid/right/.style={draw=black, fill=white},
    cuboid/top/.style={draw=black, fill=white}
  }
  \tikzcuboid{
    dimx=3, dimy=3, dimz=3,
    shiftx=50, shifty=40
  };
  \colorcubiefront[shiftx=50,shifty=40,dimx=3,dimy=3,dimz=3]{2}{3}{B}
  \colorcubiefront[shiftx=50,shifty=40,dimx=3,dimy=3,dimz=3]{3}{2}{B}
  \colorcubiefront[shiftx=50,shifty=40,dimx=3,dimy=3,dimz=3]{1}{2}{B}
  \colorcubiefront[shiftx=50,shifty=40,dimx=3,dimy=3,dimz=3]{2}{1}{B}
  \colorcubietop[shiftx=50,shifty=40,dimx=3,dimy=3,dimz=3]{2}{3}{G}
  \colorcubietop[shiftx=50,shifty=40,dimx=3,dimy=3,dimz=3]{3}{2}{G}
  \colorcubietop[shiftx=50,shifty=40,dimx=3,dimy=3,dimz=3]{1}{2}{G}
  \colorcubietop[shiftx=50,shifty=40,dimx=3,dimy=3,dimz=3]{2}{1}{G}
  \colorcubieright[shiftx=50,shifty=40,dimx=3,dimy=3,dimz=3]{2}{3}{R}
  \colorcubieright[shiftx=50,shifty=40,dimx=3,dimy=3,dimz=3]{3}{2}{R}
  \colorcubieright[shiftx=50,shifty=40,dimx=3,dimy=3,dimz=3]{1}{2}{R}
  \colorcubieright[shiftx=50,shifty=40,dimx=3,dimy=3,dimz=3]{2}{1}{R}
\end{tikzpicture}   }}. 
\]

\begin{remark}
    \label{rem:R3}
The abelianization of $\Rub_3$ is $\Cyc{2}$ and its center is
generated by the nontrivial order-two-element
$ ((0, \ldots,0) , \id ),(1, \ldots,1) , \id ))$, 
corresponding to the \emph{superflip} (i.e., the cube with all edge cubies rotated).
\end{remark}

\subsection{The Revenge Cube}

For the structure of the group $\Rub_4$ we refer to~\cite{MR795248}. Let us briefly describe this group.
The cube has corner cubies (with three possible orientations). The number of center cubies is $6 \times 4 = 24$
while the number of edges is
$12\times 2=24 ,$ too.

Positions and orientations of corners can be  described by a permutation $ \sigma\in\Sym_{8} $ and a vector $ x=(x_1, \dots ,x_8)\in (\Cyc{3})^{8} $, respectively.

A configuration for the Revenge Cube is  a tuple 
$$(x, \sigma, y, \tau, \rho_{c}, \rho_e)\in (\Cyc{3})^8\times\Sym_{8}\times\Sym_{24}\times \Sym_{24}.$$

From the analysis of \cite{MR3660738} and \cite{MR795248} one finds that all valid configurations
are determined by tuples of the form 
\[
(x, \sigma, \rho_{c},  \rho_e)\in(\Cyc3)^8\times\Sym_8\times\Sym_{24}\times\Sym_{24}
\]
satisfying the  
following conditions: 
\begin{enumerate}
\item $\sum_{i=1}^8 x_{i}=0$, 
\item $\sign (\sigma)=\sign (\rho_{c}) $.
\end{enumerate}
There is no condition 
on $\rho_e$.
Note that, even though every $\Sym_{24}$ permutation is allowed for edge cubies, their orientation is determined by their position.



In summary, the group of the Revenge Cube is 
\begin{equation}
\label{eq:revenge}
\mathcal{R}_4\simeq(\Cyc{3}\wr\Sym_8)^\circ \,\times_{\sgn}\Sym_{24}\times\Sym_{24},
\end{equation}
where the fiber product is 
given by the following maps 
\begin{gather*}
    \Sym_{24}\to\{ \pm 1\},
    \quad 
    \rho_c \mapsto \sign (\rho_c),\\
    (\Cyc{3}\wr\Sym_8)^\circ \to\{ \pm 1\},
    \quad (x,\sigma) \mapsto \sign (\sigma).
\end{gather*}

To understand this group, we refer to the 
following picture:
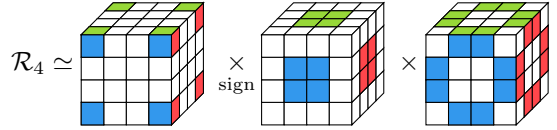
\[
\Rub_4\simeq
  \vcenter{\hbox{\begin{tikzpicture}[scale=0.3]
  \tikzset{
    cuboid/front/.style={draw=black, fill=white},
    cuboid/right/.style={draw=black, fill=white},
    cuboid/top/.style={draw=black, fill=white}
  }
  \tikzcuboid{
    dimx=4, dimy=4, dimz=4,
    shiftx=50, shifty=40
  };
  \colorcubiefront[shiftx=50,shifty=40,dimx=4,dimy=4,dimz=4]{4}{4}{B}
  \colorcubiefront[shiftx=50,shifty=40,dimx=4,dimy=4,dimz=4]{1}{1}{B}
  \colorcubiefront[shiftx=50,shifty=40,dimx=4,dimy=4,dimz=4]{1}{4}{B}
  \colorcubiefront[shiftx=50,shifty=40,dimx=4,dimy=4,dimz=4]{4}{1}{B}

  \colorcubieright[shiftx=50,shifty=40,dimx=4,dimy=4,dimz=4]{1}{1}{R}
  \colorcubieright[shiftx=50,shifty=40,dimx=4,dimy=4,dimz=4]{1}{4}{R}
  \colorcubieright[shiftx=50,shifty=40,dimx=4,dimy=4,dimz=4]{4}{1}{R}
  \colorcubieright[shiftx=50,shifty=40,dimx=4,dimy=4,dimz=4]{4}{4}{R}
  
  \colorcubietop[shiftx=50,shifty=40,dimx=4,dimy=4,dimz=4]{1}{1}{G}
  \colorcubietop[shiftx=50,shifty=40,dimx=4,dimy=4,dimz=4]{1}{4}{G}
  \colorcubietop[shiftx=50,shifty=40,dimx=4,dimy=4,dimz=4]{4}{1}{G}
  \colorcubietop[shiftx=50,shifty=40,dimx=4,dimy=4,dimz=4]{4}{4}{G} 

\end{tikzpicture}   }}
\Xsub{\sign}
  \vcenter{\hbox{\begin{tikzpicture}[scale=0.3]
  \tikzset{
    cuboid/front/.style={draw=black, fill=white},
    cuboid/right/.style={draw=black, fill=white},
    cuboid/top/.style={draw=black, fill=white}
  }
  \tikzcuboid{
    dimx=4, dimy=4, dimz=4,
    shiftx=50, shifty=40
  };
  \colorcubiefront[shiftx=50,shifty=40,dimx=4,dimy=4,dimz=4]{3}{3}{B}
  \colorcubiefront[shiftx=50,shifty=40,dimx=4,dimy=4,dimz=4]{2}{2}{B}
  \colorcubiefront[shiftx=50,shifty=40,dimx=4,dimy=4,dimz=4]{2}{3}{B}
  \colorcubiefront[shiftx=50,shifty=40,dimx=4,dimy=4,dimz=4]{3}{2}{B}

  \colorcubieright[shiftx=50,shifty=40,dimx=4,dimy=4,dimz=4]{2}{2}{R}
  \colorcubieright[shiftx=50,shifty=40,dimx=4,dimy=4,dimz=4]{2}{3}{R}
  \colorcubieright[shiftx=50,shifty=40,dimx=4,dimy=4,dimz=4]{3}{2}{R}
  \colorcubieright[shiftx=50,shifty=40,dimx=4,dimy=4,dimz=4]{3}{3}{R}
  
  \colorcubietop[shiftx=50,shifty=40,dimx=4,dimy=4,dimz=4]{2}{2}{G}
  \colorcubietop[shiftx=50,shifty=40,dimx=4,dimy=4,dimz=4]{2}{3}{G}
  \colorcubietop[shiftx=50,shifty=40,dimx=4,dimy=4,dimz=4]{3}{2}{G}
  \colorcubietop[shiftx=50,shifty=40,dimx=4,dimy=4,dimz=4]{3}{3}{G} 

\end{tikzpicture}   }}
  \times
  \vcenter{\hbox{\begin{tikzpicture}[scale=0.3]
  \tikzset{
    cuboid/front/.style={draw=black, fill=white},
    cuboid/right/.style={draw=black, fill=white},
    cuboid/top/.style={draw=black, fill=white}
  }
  \tikzcuboid{
    dimx=4, dimy=4, dimz=4,
    shiftx=50, shifty=40
  };
  \colorcubiefront[shiftx=50,shifty=40,dimx=4,dimy=4,dimz=4]{4}{3}{B}
  \colorcubiefront[shiftx=50,shifty=40,dimx=4,dimy=4,dimz=4]{4}{2}{B}
  \colorcubiefront[shiftx=50,shifty=40,dimx=4,dimy=4,dimz=4]{2}{4}{B}
  \colorcubiefront[shiftx=50,shifty=40,dimx=4,dimy=4,dimz=4]{3}{4}{B}

  \colorcubieright[shiftx=50,shifty=40,dimx=4,dimy=4,dimz=4]{4}{2}{R}
  \colorcubieright[shiftx=50,shifty=40,dimx=4,dimy=4,dimz=4]{4}{3}{R}
  \colorcubieright[shiftx=50,shifty=40,dimx=4,dimy=4,dimz=4]{3}{4}{R}
  \colorcubieright[shiftx=50,shifty=40,dimx=4,dimy=4,dimz=4]{2}{4}{R}
  
  \colorcubietop[shiftx=50,shifty=40,dimx=4,dimy=4,dimz=4]{4}{2}{G}
  \colorcubietop[shiftx=50,shifty=40,dimx=4,dimy=4,dimz=4]{4}{3}{G}
  \colorcubietop[shiftx=50,shifty=40,dimx=4,dimy=4,dimz=4]{2}{4}{G}
  \colorcubietop[shiftx=50,shifty=40,dimx=4,dimy=4,dimz=4]{3}{4}{G}

  \colorcubiefront[shiftx=50,shifty=40,dimx=4,dimy=4,dimz=4]{1}{3}{B}
  \colorcubiefront[shiftx=50,shifty=40,dimx=4,dimy=4,dimz=4]{1}{2}{B}
  \colorcubiefront[shiftx=50,shifty=40,dimx=4,dimy=4,dimz=4]{2}{1}{B}
  \colorcubiefront[shiftx=50,shifty=40,dimx=4,dimy=4,dimz=4]{3}{1}{B}

  \colorcubieright[shiftx=50,shifty=40,dimx=4,dimy=4,dimz=4]{1}{2}{R}
  \colorcubieright[shiftx=50,shifty=40,dimx=4,dimy=4,dimz=4]{1}{3}{R}
  \colorcubieright[shiftx=50,shifty=40,dimx=4,dimy=4,dimz=4]{3}{1}{R}
  \colorcubieright[shiftx=50,shifty=40,dimx=4,dimy=4,dimz=4]{2}{1}{R}
  
  \colorcubietop[shiftx=50,shifty=40,dimx=4,dimy=4,dimz=4]{1}{2}{G}
  \colorcubietop[shiftx=50,shifty=40,dimx=4,dimy=4,dimz=4]{1}{3}{G}
  \colorcubietop[shiftx=50,shifty=40,dimx=4,dimy=4,dimz=4]{2}{1}{G}
  \colorcubietop[shiftx=50,shifty=40,dimx=4,dimy=4,dimz=4]{3}{1}{G}
  
\end{tikzpicture}   }}.
\]

The initial configuration is the tuple where all permutations are the respective identity, and any vector's component equals zero. A configuration is valid if and only if it belongs to the orbit of the initial configuration, that is, it can be obtained by the application of a finite sequence of moves from the initial configuration.

\subsection{The Professor's Cube}

The cube has corner cubies (with three possible orientations). The number of interior cubies is $6 \times 9 = 54$, of which six (the center of each face) are fixed. 

Positions and orientations of corners can be described by a permutation $ \sigma\in\Sym_{8} $ and a vector 
$x=(x_1, \dots ,x_8)\in (\Cyc{3})^{8} $, respectively.

A configuration for the Professor's Cube is a tuple 
$$(x, \sigma_c, y, \sigma_e, \tau, \rho_{c}, \rho_e)\in (\Cyc{3})^8\times\Sym_{8}\times
(\Cyc{2})^{12}\times\Sym_{12}
\times\Sym_{24}\times \Sym_{24}\times \Sym_{24}.$$
This configuration is valid
if and only if the following conditions hold:
\begin{enumerate}
    \item $\sum_{i=1}^8 x_i=0$.
    \item $\sum_{i=1}^{12} y_i=0$.
    \item $\sgn(\sigma_c)=\sgn(\sigma_e)
    =\sgn(\tau)$.
    \item $\sgn(\tau)=\sgn(\rho_c)\sgn(\rho_e)$.
\end{enumerate}

Following the analysis of \cite{arXiv:2112.08602} and
\cite{bonzio_loi_peruzzi_2018}, we obtain that
the group $\Rub_5$ is isomorphic~to 
\[
\mathcal{R}_5=\
\left(\mathcal{R}_3
\times_{\sgn} \Sym_{24}
\right)
\times_{\sgn}
(\Sym_{24}\ \times \Sym_{24} ),
\]
where the fiber products are 
given by the following maps
\begin{align*}
    &\Rub_3\to\{\pm1\}, 
    &&(x,\sigma_c,y,\sigma_e)
    \mapsto
    \sgn(\sigma_c)
    ,\\
    &\Sym_{24}\to\{\pm1\},
    &&\tau \mapsto\sgn(\tau ),
\shortintertext{(recall that $\sgn(\sigma_e)= \sgn(\sigma_c)$) and}
    &\mathcal{R}_3\times_{\sgn} \Sym_{24}\mapsto\{\pm1\},
    &&(x,\sigma_c,y,\sigma_e,\tau)\mapsto 
    \sgn(\tau) 
    = \sgn(\sigma_e), 
    \\
    &\Sym_{24}\times\Sym_{24}\to\{\pm1\},
    &&(\rho_c,\rho_e)\mapsto\sgn(\rho_c),
\end{align*}
respectively (note that $\sgn(\tau)= \sgn(\sigma_e)$). 

To understand how the group $\Rub_3$ and the copies of the symmetric group $\Sym_{24}$ arise in this description of $\Rub_5$, we refer to the following picture: 
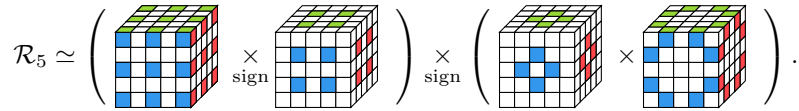
\[
\Rub_5\simeq\left(\
  \vcenter{\hbox{\begin{tikzpicture}[scale=0.2]
  \tikzset{
    cuboid/front/.style={draw=black, fill=white},
    cuboid/right/.style={draw=black, fill=white},
    cuboid/top/.style={draw=black, fill=white}
  }
  \tikzcuboid{
    dimx=5, dimy=5, dimz=5,
    shiftx=50, shifty=40
  };
  \colorcubiefront[shiftx=50,shifty=40,dimx=5,dimy=5,dimz=5]{5}{1}{B}
  \colorcubiefront[shiftx=50,shifty=40,dimx=5,dimy=5,dimz=5]{3}{3}{B}
  \colorcubiefront[shiftx=50,shifty=40,dimx=5,dimy=5,dimz=5]{1}{1}{B}
  \colorcubiefront[shiftx=50,shifty=40,dimx=5,dimy=5,dimz=5]{1}{5}{B}
  \colorcubiefront[shiftx=50,shifty=40,dimx=5,dimy=5,dimz=5]{5}{5}{B}
  
  \colorcubietop[shiftx=50,shifty=40,dimx=5,dimy=5,dimz=5]{1}{1}{G}
  \colorcubietop[shiftx=50,shifty=40,dimx=5,dimy=5,dimz=5]{1}{5}{G}
  \colorcubietop[shiftx=50,shifty=40,dimx=5,dimy=5,dimz=5]{5}{5}{G}
  \colorcubietop[shiftx=50,shifty=40,dimx=5,dimy=5,dimz=5]{5}{1}{G}
  \colorcubietop[shiftx=50,shifty=40,dimx=5,dimy=5,dimz=5]{3}{3}{G}

  \colorcubieright[shiftx=50,shifty=40,dimx=5,dimy=5,dimz=5]{1}{1}{R}
  \colorcubieright[shiftx=50,shifty=40,dimx=5,dimy=5,dimz=5]{1}{5}{R}
  \colorcubieright[shiftx=50,shifty=40,dimx=5,dimy=5,dimz=5]{5}{5}{R}
  \colorcubieright[shiftx=50,shifty=40,dimx=5,dimy=5,dimz=5]{5}{1}{R}
  \colorcubieright[shiftx=50,shifty=40,dimx=5,dimy=5,dimz=5]{3}{3}{R}

  \colorcubieright[shiftx=50,shifty=40,dimx=5,dimy=5,dimz=5]{1}{1}{R}
  \colorcubietop[shiftx=50,shifty=40,dimx=5,dimy=5,dimz=5]{1}{1}{G}

  \colorcubiefront[shiftx=50,shifty=40,dimx=5,dimy=5,dimz=5]{5}{3}{B}
  \colorcubiefront[shiftx=50,shifty=40,dimx=5,dimy=5,dimz=5]{3}{5}{B}
  \colorcubiefront[shiftx=50,shifty=40,dimx=5,dimy=5,dimz=5]{1}{3}{B}
  \colorcubiefront[shiftx=50,shifty=40,dimx=5,dimy=5,dimz=5]{3}{1}{B}

  \colorcubieright[shiftx=50,shifty=40,dimx=5,dimy=5,dimz=5]{5}{3}{R}
  \colorcubieright[shiftx=50,shifty=40,dimx=5,dimy=5,dimz=5]{3}{5}{R}
  \colorcubieright[shiftx=50,shifty=40,dimx=5,dimy=5,dimz=5]{1}{3}{R}
  \colorcubieright[shiftx=50,shifty=40,dimx=5,dimy=5,dimz=5]{3}{1}{R}
  
  \colorcubietop[shiftx=50,shifty=40,dimx=5,dimy=5,dimz=5]{5}{3}{G}
  \colorcubietop[shiftx=50,shifty=40,dimx=5,dimy=5,dimz=5]{3}{5}{G}
  \colorcubietop[shiftx=50,shifty=40,dimx=5,dimy=5,dimz=5]{1}{3}{G}
  \colorcubietop[shiftx=50,shifty=40,dimx=5,dimy=5,dimz=5]{3}{1}{G}
\end{tikzpicture}   }}
  \Xsub{\sign}
  \vcenter{\hbox{\begin{tikzpicture}[scale=0.2]
  \tikzset{
    cuboid/front/.style={draw=black, fill=white},
    cuboid/right/.style={draw=black, fill=white},
    cuboid/top/.style={draw=black, fill=white}
  }
  \tikzcuboid{
    dimx=5, dimy=5, dimz=5,
    shiftx=50, shifty=40
  };
  \colorcubiefront[shiftx=50,shifty=40,dimx=5,dimy=5,dimz=5]{4}{4}{B}
  \colorcubiefront[shiftx=50,shifty=40,dimx=5,dimy=5,dimz=5]{2}{4}{B}
  \colorcubiefront[shiftx=50,shifty=40,dimx=5,dimy=5,dimz=5]{4}{2}{B}
  \colorcubiefront[shiftx=50,shifty=40,dimx=5,dimy=5,dimz=5]{2}{2}{B}
  
  \colorcubietop[shiftx=50,shifty=40,dimx=5,dimy=5,dimz=5]{4}{4}{G}
  \colorcubietop[shiftx=50,shifty=40,dimx=5,dimy=5,dimz=5]{2}{4}{G}
  \colorcubietop[shiftx=50,shifty=40,dimx=5,dimy=5,dimz=5]{4}{2}{G}
  \colorcubietop[shiftx=50,shifty=40,dimx=5,dimy=5,dimz=5]{2}{2}{G}
  
  \colorcubieright[shiftx=50,shifty=40,dimx=5,dimy=5,dimz=5]{4}{4}{R}
  \colorcubieright[shiftx=50,shifty=40,dimx=5,dimy=5,dimz=5]{2}{4}{R}
  \colorcubieright[shiftx=50,shifty=40,dimx=5,dimy=5,dimz=5]{4}{2}{R}
  \colorcubieright[shiftx=50,shifty=40,dimx=5,dimy=5,dimz=5]{2}{2}{R}
\end{tikzpicture}}}\ 
\right)
\Xsub{\sign}
\left(\
  \vcenter{\hbox{\begin{tikzpicture}[scale=0.2]
  \tikzset{
    cuboid/front/.style={draw=black, fill=white},
    cuboid/right/.style={draw=black, fill=white},
    cuboid/top/.style={draw=black, fill=white}
  }
  \tikzcuboid{
    dimx=5, dimy=5, dimz=5,
    shiftx=50, shifty=40
  };
  \colorcubiefront[shiftx=50,shifty=40,dimx=5,dimy=5,dimz=5]{4}{3}{B}
  \colorcubiefront[shiftx=50,shifty=40,dimx=5,dimy=5,dimz=5]{3}{4}{B}
  \colorcubiefront[shiftx=50,shifty=40,dimx=5,dimy=5,dimz=5]{2}{3}{B}
  \colorcubiefront[shiftx=50,shifty=40,dimx=5,dimy=5,dimz=5]{3}{2}{B}

  \colorcubieright[shiftx=50,shifty=40,dimx=5,dimy=5,dimz=5]{4}{3}{R}
  \colorcubieright[shiftx=50,shifty=40,dimx=5,dimy=5,dimz=5]{3}{4}{R}
  \colorcubieright[shiftx=50,shifty=40,dimx=5,dimy=5,dimz=5]{2}{3}{R}
  \colorcubieright[shiftx=50,shifty=40,dimx=5,dimy=5,dimz=5]{3}{2}{R}
  
  \colorcubietop[shiftx=50,shifty=40,dimx=5,dimy=5,dimz=5]{4}{3}{G}
  \colorcubietop[shiftx=50,shifty=40,dimx=5,dimy=5,dimz=5]{3}{4}{G}
  \colorcubietop[shiftx=50,shifty=40,dimx=5,dimy=5,dimz=5]{2}{3}{G}
  \colorcubietop[shiftx=50,shifty=40,dimx=5,dimy=5,dimz=5]{3}{2}{G}
  
\end{tikzpicture}}}
  \times
  \vcenter{\hbox{\begin{tikzpicture}[scale=0.2]
  \tikzset{
    cuboid/front/.style={draw=black, fill=white},
    cuboid/right/.style={draw=black, fill=white},
    cuboid/top/.style={draw=black, fill=white}
  }
  \tikzcuboid{
    dimx=5, dimy=5, dimz=5,
    shiftx=50, shifty=40
  };
  \colorcubiefront[shiftx=50,shifty=40,dimx=5,dimy=5,dimz=5]{4}{1}{B}
  \colorcubiefront[shiftx=50,shifty=40,dimx=5,dimy=5,dimz=5]{2}{1}{B}
  \colorcubiefront[shiftx=50,shifty=40,dimx=5,dimy=5,dimz=5]{1}{4}{B}
  \colorcubiefront[shiftx=50,shifty=40,dimx=5,dimy=5,dimz=5]{1}{2}{B}
  \colorcubiefront[shiftx=50,shifty=40,dimx=5,dimy=5,dimz=5]{5}{2}{B}
  \colorcubiefront[shiftx=50,shifty=40,dimx=5,dimy=5,dimz=5]{5}{4}{B}
  \colorcubiefront[shiftx=50,shifty=40,dimx=5,dimy=5,dimz=5]{2}{5}{B}
  \colorcubiefront[shiftx=50,shifty=40,dimx=5,dimy=5,dimz=5]{4}{5}{B}

  \colorcubietop[shiftx=50,shifty=40,dimx=5,dimy=5,dimz=5]{4}{1}{G}
  \colorcubietop[shiftx=50,shifty=40,dimx=5,dimy=5,dimz=5]{2}{1}{G}
  \colorcubietop[shiftx=50,shifty=40,dimx=5,dimy=5,dimz=5]{1}{4}{G}
  \colorcubietop[shiftx=50,shifty=40,dimx=5,dimy=5,dimz=5]{1}{2}{G}
  \colorcubietop[shiftx=50,shifty=40,dimx=5,dimy=5,dimz=5]{5}{2}{G}
  \colorcubietop[shiftx=50,shifty=40,dimx=5,dimy=5,dimz=5]{5}{4}{G}
  \colorcubietop[shiftx=50,shifty=40,dimx=5,dimy=5,dimz=5]{2}{5}{G}
  \colorcubietop[shiftx=50,shifty=40,dimx=5,dimy=5,dimz=5]{4}{5}{G}
  
  \colorcubieright[shiftx=50,shifty=40,dimx=5,dimy=5,dimz=5]{4}{1}{R}
  \colorcubieright[shiftx=50,shifty=40,dimx=5,dimy=5,dimz=5]{2}{1}{R}
  \colorcubieright[shiftx=50,shifty=40,dimx=5,dimy=5,dimz=5]{1}{4}{R}
  \colorcubieright[shiftx=50,shifty=40,dimx=5,dimy=5,dimz=5]{1}{2}{R}
  \colorcubieright[shiftx=50,shifty=40,dimx=5,dimy=5,dimz=5]{5}{2}{R}
  \colorcubieright[shiftx=50,shifty=40,dimx=5,dimy=5,dimz=5]{5}{4}{R}
  \colorcubieright[shiftx=50,shifty=40,dimx=5,dimy=5,dimz=5]{2}{5}{R}
  \colorcubieright[shiftx=50,shifty=40,dimx=5,dimy=5,dimz=5]{4}{5}{R}

\end{tikzpicture}}}\
\right).
\]

\section{The proofs}
\label{section:proofs}

The proofs of Theorems~\ref{thm:n=4} and~\ref{thm:n=5} rely heavily on the structure of the group $\Rub_3$. 

%
%
%
%

\begin{proof}[Proof of Theorem~\ref{thm:n=4}]



We need two linearly disjoint extensions
with Galois groups
$  \S_{24} \times_{\sgn} (\Cyc{3} 
\wr \S_{8})^\circ $
and $\S_{24}$, respectively. 
For the first such extension, let $f$ as in Theorem~\ref{thm:Rubik} with
parameters $a$ and $b$ as in~\eqref{eq:Rubik_parameters}, 
and 
\[
g=X^{24} - t(X+1)
\]
be such that 
$\disc{g} \equiv
7 \cdot 1437417619559484462138047
\bmod 
(\Q^\times)^2$. These polynomials come from Theorem~\ref{thm:Rubik} 
and the Nart--Vila theorem~\cite{zbMATH01458933,zbMATH01458934}, respectively.  



Setting
$$
t = \frac{24^{24}}{23^{23} } 
\frac{-1}
{(23\cdot 7 \cdot 1437417619559484462138047
s^2 +1 )}
$$
we get
\[
\disc{g} = \frac{2^{1728}3^{576}s^2 (s^2 + 1/231424236749076998404225567)^{-24}}{7^{23}23^{552} 
1437417619559484462138047^{23}}
\]
which is congruent to 
$ 7 \cdot 1437417619559484462138047 $ 
in 
$\Q^\times/(\Q^\times)^2$ for any specialization $s\in \Q^\times.$

By the condition on the discriminants, 
$$
\Galf{fg}{\Q} \subseteq \S_{24} \times_{\sgn} ( \Cyc{3} \wr \S_8)^\circ
$$
having as quotients
$$
\Galf{f}{\Q} \simeq  ( \Cyc{3} \wr \S_8)^\circ,\quad 
\Galf{g}{\Q} \simeq \S_{24}.
$$

The Jordan--H\"older decomposition of the Galois group 
$\Galf{fg}{\Q}$ has (at least) components
$\A_{24}, \A_{8}$ and $\Cyc2$ with multiplicity one 
and $\Cyc{3}$ with multiplicity seven.

By order considerations, one gets the isomorphism
$$
\Galf{fg}{\Q} \simeq \S_{24} \times_{\sgn} ( \Cyc{3} \wr \S_8)^\circ.
$$

The same argument from \cite[Remark 1]{arXiv:2411.11566}
using the results from Nart--Vila \cite{zbMATH01458934,zbMATH01458933}
proves that 
$\Galf{h}{\Q} \simeq \Sym_{24}.$ 

To get the polynomials in the statement we take $s=1$ and   
the claim follows. 
\end{proof}


\begin{proof}[Proof of Theorem~\ref{thm:n=5}]
Recall that $f$ and $g$ generate an $\Rub_3$-field extension over $\Q.$
Its unique quadratic subextension is 
$\Q[\sqrt{
7 \cdot 1437417619559484462138047
}]$
(by Remark~\ref{rem:R3}).

Let us now take care of the other three $\S_{24}$ factors.
We need three $\S_{24}$-polynomials 
$h_i = X^{24}-u_i(X+1)$ for $i=1,2,3$ 
such that 
\begin{align}
\disc{h_1} &=
(-1)  (23^{23} u_1 + 24^{24})  u_1^{23}
\nonumber \\ &\equiv 
7 \cdot 1437417619559484462138047
\bmod 
(\Q^\times)^2,
\label{e:disch1}\\
\disc{h_2}\disc{h_3} &=  
(23^{23} u_2 + 24^{24})  u_2^{23} (23^{23} u_3 + 24^{24})  u_3^{23}
\nonumber \\
&\equiv 
7 \cdot 1437417619559484462138047
\bmod 
(\Q^\times)^2.\label{e:disch23}
\end{align}

Condition~\eqref{e:disch1} on $u_1$ takes care of the
fiber product along ${\varphi_{\sigma},\varphi_{\rho_c}}$, and
Condition~\eqref{e:disch23} 
on $u_2$ and $u_3$ of the 
fiber product with respect to ${\varphi_{\sigma}, \varphi_{\tau \rho_e}}$.



Simplifying and replacing 
by $u_i = v_i\frac{24^{24} }{ 23^{22} }$
for $i\in\{2,3\}$, 
\begin{align*}
\disc{h_2}\disc{h_3} &\equiv
(23^{23} u_2 + 24^{24})   (23^{23} u_3 + 24^{24})  u_3 u_2\bmod 
(\Q^\times)^2\\
&\equiv (23 v_2 + 1)   (23 v_3 +  1)  v_3 v_2\bmod 
(\Q^\times)^2\\
&\equiv 
7 \cdot 1437417619559484462138047
\bmod 
(\Q^\times)^2.
\end{align*}
One possible way to achieve these conditions is to consider
\begin{align*}
u_1 &= \frac{24^{24}}{23^{23} } 
\frac{-1}
{(23\cdot 7 \cdot 1437417619559484462138047
v_1^2 +1 )},\\
v_1&=1,
\end{align*}
to guarantee Condition~\eqref{e:disch1}, and
\begin{align*}
v_2 &= 7 \cdot 1437417619559484462138047,\\
v_3 &= z\cdot w^2, \\
\shortintertext{with}
w &= \frac{2r}{23 z -r^2 },\quad 
    z = 14464014796817312400264098, \quad 
    r=1, 
\end{align*}
for Condition~\eqref{e:disch23}.


The conditions on the discriminants guarantee that the Galois group of $h_1h_2h_3$ is a subgroup of  
\[
\Sym_{24}
\times_{\sgn}
(\Sym_{24}\ \times \Sym_{24} ). 
\]
Magma shows that
$\Galf{h_i}{\Q}\simeq \S_{24}$ for all $i \in\{1,2,3\}$. Moreover, it verifies that the polynomial $h_1h_2h_3$ has Galois group of order $(24!)^3/2$. On our desktop computer, this computation requires approximately 18 hours.
\end{proof}



\appendix 

\section{Magma code for verifying theorems}

The proofs of Theorems~\ref{thm:n=4} and~\ref{thm:n=5} are entirely theoretical and do not require any computer assistance. However, the statements can be verified computationally in a straightforward manner using Magma. As an example, we present a script that verifies Theorem~\ref{thm:n=4} for a particular choice of parameters. 

\begin{lstlisting}
> Q<x>:=PolynomialAlgebra(Rationals());
> t := 2^67*3^24/(23^23*31*281*1201*70529*9801219477271);
> g := x^24-t*(x+1);
> a := 2139;
> b := -6489;
> f8 := x^8-a*x-b;
> f := Q!((x^2-x)^8*Evaluate(f8, (x^3-3*x+1)/(x^2-x)));
> time H := GaloisGroup(f*g*(x^24-x-1));
Time: 666.670
> load rubik444;
Loading "/Applications/Magma/libs/pergps/rubik444"
The automorphism group of the 4 x 4 x 4 Rubik cube.
The group is represented as a permutation group of degree 72.
Its order is
2^50 * 3^29 * 5^9 * 7^7 * 11^4 * 13^2 * 17^2 * 19^2 * 23^2.
Group: G
> time flag, f := IsIsomorphic(G,H);
Time: 5.170
> flag
true
\end{lstlisting}

Magma provides a presentation of the Revenge Cube group, which is loaded on line 10. On line 17, we 
verify that the Galois group we constructed is isomorphic 
to this group. 
This line of the code is written this way solely to avoid printing an explicit (and lengthy) isomorphism between the groups.

An analogous computation applies to Theorem~\ref{thm:n=5} and is omitted.

\section{Generators for the Professor's Cube group}
\label{section:generators}

Generators for the groups $\Rub_3$ and $\Rub_4$ are well known and are even available in Magma. Generators for the group $\Rub_5$ can be obtained in the same way as the generators for the other groups: we label each cubie and record how each move permutes the labels; see Figure~\ref{fig:net5x5x5}.

Although presenting these generators requires considerable space, the fact that they are difficult to find in the literature (and, unfortunately, that those which can be found are typically incorrect) suggests that it is time to write them down once and for all. 

A generating set of $\Rub_5$ is given by the permutations  
\begin{align*}
&r_1 = (40\,88\,9\,96)(28\,76\,21\,84)(16\,64\,33\,72)(4\,52\,45\,60)(41\,5\,8\,44)\\
&\qquad\qquad (42\,17\,7\,32)(43\,29\,6\,20)(18\,19\,31\,30)(101\,127\,104\,126)\\
&\qquad\qquad (98\,128\,107\,125)(124\,136\,129\,144),\\
&r_2 = (39\,87\,10\,95)(27\,75\,22\,83)(15\,63\,34\,71)(3\,51\,46\,59)(123\,135\,130\,143),\\
&b_1 = (52\,53\,93\,44)(51\,65\,94\,32)(50\,77\,95\,20)(49\,89\,96\,8)(9\,12\,48\,45)\\
&\qquad\qquad (10\,24\,47\,33)(11\,36\,46\,21)(22\,23\,35\,34)(99\,132\,108\,129)\\
&\qquad\qquad (102\,131\,105\,130)(109\,137\,120\,128),\\
&b_2 = (54\,81\,43\,64)(66\,82\,31\,63)(78\,83\,19\,62)(90\,84\,7\,61)(138\,117\,127\,112),\\
&d_1 = (57\,60\,96\,93)(58\,72\,95\,81)(59\,84\,94\,69)(70\,71\,83\,82)(45\,89\,37\,41)\\
&\qquad\qquad (46\,90\,38\,42)(47\,91\,39\,43)(48\,92\,40\,44)(111\,144\,120\,141)\\
&\qquad\qquad (114\,143\,117\,142)(108\,119\,106\,107),\\
&d_2 = (33\,77\,25\,29)(34\,78\,26\,30)(35\,79\,27\,31)(36\,80\,28\,32)(105\,116\,103\,104),\\
&u_1 = (49\,52\,88\,85)(62\,63\,75\,74)(50\,64\,87\,73)(51\,76\,86\,61)(5\,1\,53\,9)\\
&\qquad\qquad (6\,2\,54\,10)(7\,3\,55\,11)(8\,4\,56\,12)(109\,136\,118\,133)\\
&\qquad\qquad (112\,135\,115\,134)(98\,97\,110\,99),\\
&u_2 = (17\,13\,65\,21)(18\,14\,66\,22)(19\,15\,67\,23)(20\,16\,68\,24)(101\,100\,113\,102),\\
&l_1 = (57\,48\,49\,1)(69\,36\,61\,13)(81\,24\,73\,25)(93\,12\,85\,37)(89\,53\,56\,92)\\
&\qquad\qquad (90\,65\,55\,80)(91\,77\,54\,68)(66\,67\,79\,78)(110\,140\,119\,137)\\
&\qquad\qquad (113\,139\,116\,138)(132\,133\,121\,141),\\
&l_2 = (94\,11\,86\,38)(82\,23\,74\,26)(70\,35\,62\,14)(58\,47\,50\,2)(131\,134\,122\,142),\\
&f_1 = (85\,5\,60\,92)(86\,17\,59\,80)(87\,29\,58\,68)(88\,41\,57\,56)(1\,4\,40\,37)\\
&\qquad\qquad (2\,16\,39\,25)(3\,28\,38\,13)(14\,15\,27\,26)(100\,123\,103\,122)\\
&\qquad\qquad (97\,124\,106\,121)(118\,125\,111\,140),\\
&f_2 = (73\,6\,72\,91)(74\,18\,71\,79)(75\,30\,70\,67)(76\,42\,69\,55)(115\,126\,114\,139).    
\end{align*}

\begin{figure}
\begin{tikzpicture}[scale=0.57]


\drawFaceRows{0}{0}{L}{FaceL}%
  {53,54,110,55,56}%
  {65,66,113,67,68}%
  {137,138,L,139,140}%
  {77,78,116,79,80}%
  {89,90,119,91,92};

\drawFaceRows{5}{0}{F}{FaceF}%
  {1,2,97,3,4}%
  {13,14,100,15,16}%
  {121,122,F,123,124}%
  {25,26,103,27,28}%
  {37,38,106,39,40};

\drawFaceRows{10}{0}{R}{FaceR}%
  {5,6,98,7,8}%
  {17,18,101,19,20}%
  {125,126,R,127,128}%
  {29,30,104,31,32}%
  {41,42,107,43,44};

\drawFaceRows{15}{0}{B}{FaceB}%
  {9,10,99,11,12}%
  {21,22,102,23,24}%
  {129,130,B,131,132}%
  {33,34,105,35,36}%
  {45,46,108,47,48};

\drawFaceRows{5}{5}{U}{FaceU}%
  {49,50,109,51,52}%
  {61,62,112,63,64}%
  {133,134,U,135,136}%
  {73,74,115,75,76}%
  {85,86,118,87,88};

\drawFaceRows{5}{-5}{D}{FaceD}%
  {57,58,111,59,60}%
  {69,70,114,71,72}%
  {141,142,D,143,144}%
  {81,82,117,83,84}%
  {93,94,120,95,96};

\end{tikzpicture}
\caption{Unfolded net of Professor's cube, with labels. Our labeling is obtained by extending one particular labeling (slightly different from the one mentioned before) of the Revenge Cube, assigning numbers starting from 97 to the central cross of each face.}
\label{fig:net5x5x5}
\end{figure}

A generating set for $\Rub_4$ can be obtained from the generating set for $\Rub_5$ that we have just presented: it suffices to disregard all elements whose labels are greater than $96$.

\section*{Acknowledgments}

A. Loi and M. Damele are supported by INdAM and  GNSAGA - Gruppo Nazionale per le Strutture Algebriche, Geometriche e le loro Applicazioni and by the project ProBiKi of Fondazione di Sardegna (Italy). M. Damele also  acknowledges the support and hospitality of the University of Valencia, where part of this work was conducted during a research stay..
Mereb is supported by grants CONICET PIP 11220210100220CO, and UBACyT 
20020220100065BA.
Vendramin is supported by OZR3762
of Vrije Universiteit Brussel
and FWO Senior Research Project G004124N.

\bibliographystyle{abbrv}
\bibliography{refs}

\begin{thebibliography}{10}

\bibitem{bandelow2012inside}
C.~Bandelow.
\newblock Inside {Rubik}'s cube and beyond, with 21 cartoons by {Alexander} {Maga}. {Transl}. from the {German} by {Jeannette} {Zehnder} and {Lucy} {Moser}.
\newblock Boston-{Basel}-{Stuttgart}: {Birkh{\"a}user}. 125 p. {DM} 18.80 (1982), 1982.

\bibitem{MR3660738}
S.~Bonzio, A.~Loi, and L.~Peruzzi.
\newblock The first law of cubology for the {R}ubik's revenge.
\newblock {\em Math. Slovaca}, 67(3):561--572, 2017.

\bibitem{bonzio_loi_peruzzi_2018}
S.~Bonzio, A.~Loi, and L.~Peruzzi.
\newblock On the n × n × n rubik's cube.
\newblock {\em Mathematica Slovaca}, 68, 10 2018.

\bibitem{zbMATH01077111}
W.~Bosma, J.~Cannon, and C.~Playoust.
\newblock The {Magma} algebra system. {I}: {The} user language.
\newblock {\em J. Symb. Comput.}, 24(3-4):235--265, 1997.

\bibitem{MR2599606}
D.~Joyner.
\newblock {\em Adventures in group theory}.
\newblock Johns Hopkins University Press, Baltimore, MD, second edition, 2008.
\newblock Rubik's cube, Merlin's machine, and other mathematical toys.

\bibitem{MR795248}
M.~E. Larsen.
\newblock Rubik's revenge: the group theoretical solution.
\newblock {\em Amer. Math. Monthly}, 92(6):381--390, 1985.

\bibitem{MR3822366}
G.~Malle and B.~H. Matzat.
\newblock {\em Inverse {G}alois theory}.
\newblock Springer Monographs in Mathematics. Springer, Berlin, 2018.
\newblock Second edition.

\bibitem{arXiv:2411.11566}
M.~Mereb and L.~Vendramin.
\newblock Rubik's as a {Galois}'.
\newblock Preprint, {arXiv}:2411.11566 [math.{NT}] (2024), 2024.

\bibitem{zbMATH01458933}
E.~Nart and N.~Vila.
\newblock Equations of the type {{\(X^{n}+aX^{2}+bX+c,\)}} {{\(n\)}} an odd square, with absolute galois group {{\(A_{n}\)}}.
\newblock In {\em Actas de las VI jornadas de matemáticas hispano-lusas. Part II}, pages 826--828. Santander: Univ. Santander, rev. {Univ}. {Santander} {No}. 2, {Part} 2 edition, 1979.

\bibitem{zbMATH01458934}
E.~Nart and N.~Vila.
\newblock Equations of the type {{\(X^n+aX+b\)}} with absolute galois group {{\(S_n\)}}.
\newblock In {\em Actas de las VI jornadas de matemáticas hispano-lusas. Part II}, pages 821--825. Santander: Univ. Santander, rev. {Univ}. {Santander} {No}. 2, {Part} 2 edition, 1979.

\bibitem{zbMATH06212024}
T.~Rokicki, H.~Kociemba, M.~Davidson, and J.~Dethridge.
\newblock The diameter of the {Rubik}'s {Cube} group is twenty.
\newblock {\em SIAM J. Discrete Math.}, 27(2):1082--1105, 2013.

\bibitem{arXiv:2112.08602}
D.~Salkinder.
\newblock $n {{\times}} n {{\times}} n$ {Rubik}'s {Cubes} and {God}'s {Number}.
\newblock Preprint, {arXiv}:2112.08602 [math.{CO}] (2021), 2021.

\bibitem{MR2363329}
J.-P. Serre.
\newblock {\em Topics in {G}alois theory}, volume~1 of {\em Research Notes in Mathematics}.
\newblock A K Peters, Ltd., Wellesley, MA, second edition, 2008.
\newblock With notes by Henri Darmon.

\bibitem{MR1405612}
H.~V\"{o}lklein.
\newblock {\em Groups as {G}alois groups}, volume~53 of {\em Cambridge Studies in Advanced Mathematics}.
\newblock Cambridge University Press, Cambridge, 1996.
\newblock An introduction.

\end{thebibliography}

\end{document}